\documentstyle[12pt]{article} 
\newlength{\abstractwidth} 
\renewcommand{\thefootnote}{\fnsymbol{footnote}}
\renewcommand{\thanks}[1]{\footnote{#1}} 
\newcommand{\starttext}{ \setcounter{footnote}{0}
\renewcommand{\thefootnote}{\arabic{footnote}}}

\newcommand{\be}{\begin{equation}}
\newcommand{\bea}{\begin{eqnarray}}
\newcommand{\eea}{\end{eqnarray}} \newcommand{\ee}{\end{equation}}
 
 \def\ba{\begin{eqnarray}}
\def\ea{\end{eqnarray}}

\newtheorem{theorem}{Theorem}

\def\log{\,{\rm log}\,}

  \def\o{\omega}

\def\al{\alpha}

\def\g{\gamma}

\def\e{\epsilon}

\def\o{\omega}

 \def\I{\int}    \def\ddb{{\partial\bar\partial}}

 \def\v{\vskip .1in}

\def\[{{\bf [}} \def\]{{\bf ]}}   \def\pl{\partial}

\begin{document} \starttext \baselineskip=16pt
\setcounter{footnote}{0}
\begin{center}
{\Large\bf THE MOSER-TRUDINGER INEQUALITY}\\
{\Large\bf ON K\"AHLER-EINSTEIN MANIFOLDS}
\footnote{Research supported in
part by National Science Foundation grants DMS-02-45371,
DMS-05-14003 and DMS-05-04285.}

\bigskip
\bigskip

\begin{center}
\begin{tabular}{ccc}
{\bf D.H. Phong}  & & {\bf Jian Song} \\
Department of Mathematics & &Department of Mathematics \\
Columbia University & & Johns Hopkins University \\
New York, NY 10027 & & Baltimore, MD 21218 \\ && \\
{\bf Jacob Sturm}  & & {\bf Ben Weinkove\footnotemark} \\
Department of Mathematics & & Department of Mathematics \\
Rutgers University & & Imperial College \\
Newark, NJ 07102 & & London, SW7 2AZ, U.K.
\end{tabular}
\end{center}
\footnotetext{Currently on leave from Harvard
University and supported by a Royal Society Research Assistantship.}
\end{center}

\v
\v
\v

{\bf ABSTRACT:} \ We prove the conjecture of Tian on the strong form of the
Moser-Trudinger inequality for K\"ahler-Einstein manifolds with positive
first Chern class, when there are no holomorphic vector fields, and, more
generally, when the setting is invariant under a maximal compact subgroup
of the automorphism group.

\bigskip

\baselineskip=13pt
\setcounter{equation}{0}
\setcounter{footnote}{0}

\section{Introduction}
\setcounter{equation}{0}

The existence of canonical metrics in differential geometry
is closely related to the properness of corresponding energy
functionals. A basic example is given by Hermitian-Einstein metrics on holomorphic vector bundles.  By the Donaldson-Uhlenbeck-Yau theorem \cite{D1}, \cite{UY}
existence is equivalent
to the stability of the bundle.  But existence is also equivalent
  to a certain non-linear bound from below
for the Donaldson energy functional \cite{D2}.  Another example
is provided by K\"ahler-Einstein metrics.  Such a metric always
exists on a manifold $X$ with $c_1(X)<0$ \cite{Y1}, \cite{A1} and $c_1(X)=0$
\cite{Y1}.  If $c_1(X)>0$ then existence is conjectured by Yau to be equivalent
to stability in the sense of geometric invariant theory \cite{Y2}.  Relating
this to energy functionals, Tian showed in
 \cite{T97}
 that  a K\"ahler-Einstein metric  $\o_{KE}$ exists in this case if and only if the
 functional $F_{\o_{KE}}(\phi)$ is proper.   More precisely,
let $(X,\o)$ be a compact K\"ahler manifold of dimension $n$
and volume $V$ with $c_1(X)>0$, and define the functionals $F_{\o}(\phi)$
and $J_{\o}(\phi)$ by
\setlength\arraycolsep{2pt}
\bea
F_{\o}(\phi)
&=&
J_{\o}(\phi)-{1\over V}\int_X\phi\,\o^n
-\log ({1\over V}\int_X e^{h_{\o}-\phi}\o^n),
\nonumber\\
J_{\o}(\phi)
&=&
{\sqrt{-1} \over 2V}\sum_{j=0}^{n-1}\left( {n-j\over n+1} \right)  \int_X \pl\phi\wedge\bar\pl\phi
\wedge \o^{n-1-j}\wedge\o_\phi^j,
\eea
for all $\phi$ in $
{\cal P}(X,\o)=\{\phi \in C^{\infty}(X) \ | \ \o_\phi=\o+{\sqrt{-1}\over 2}\ddb\phi>0\}$.
Here $h_{\o}$ is the potential for the Ricci curvature of the
K\"ahler form $\o$ defined by
\bea
\label{hnorm}
\textrm{Ric}(\o)-\o={\sqrt{-1}\over 2}\ddb h_\o,
\qquad
\int_X e^{h_\o}\o^n=V.
\eea
Then it is shown in \cite{T97} that, under the assumption that $X$
admits no non-trivial holomorphic vector fields, the existence of
a K\"ahler-Einstein metric
$\o_{KE}$ is equivalent to the existence of a
K\"ahler metric
$\o\in c_1(X)$,  a strictly positive exponent $\g$, and  corresponding
constants $A_\g>0$, $B_\g$, so that
\bea
\label{MT}
F_{\o}(\phi)\ \geq \ A_\g\,J_{\o}(\phi)^\g\,-\,B_\g,
\eea
for all $\phi$ in ${\cal P}(X, \o)$. In view of the cocycle identity,
$F_{\o}(\phi_1+\phi_2)=F_{\o}(\phi_1)+F_{\o_{\phi_1}}(\phi_2)$ and properties
of the functional $J_{\o}$ (see section 2),
the inequality (\ref{MT}) holds for every $\o\in c_1(X)$ if and only
if it holds for some $\o\in c_1(X)$.

\v

Inequalities such as (\ref{MT})
are sometimes referred to as generalized Moser-Trudinger inequalities,
by analogy with similar identities on the sphere $S^2$.  Indeed, the Moser-Trudinger
inequality \cite{Mo}, \cite{Tr}, \cite{On} for functions $u \in L^2_1(S^2)$,
\bea \label{aubin}
\frac{1}{V} \int_{S^2} e^{-u} \o_{KE} \le C e^{\frac{\nu}{4V} \int_{S^2} \sqrt{-1}
\partial u \wedge \overline{\partial} u - \frac{1}{V} \int_{S^2} u \, \o_{KE}},
\eea
with the optimal constants $C=\nu=1$ and $\omega_{KE}$ the round metric on $S^2$,  implies $F_{\o_{KE}}(\phi)
\ge 0$ for $\phi \in {\cal P}(S^2, \o_{KE})$.  Aubin \cite{A2} showed that
if $u$ satisfies
\bea \label{orthogonal}
\int_{S^2} \lambda e^{-u}
\o_{KE} =0, \quad \textrm{for all } \lambda \in \Lambda,
\eea
where $\Lambda$ is the space of eigenfunctions of the Laplacian $\Delta_{\o_{KE}}$
with eigenvalue 1, then
 the constant $\nu$ can be reduced to $1/2 + \epsilon$ for any $\epsilon>0$, at the expense of increasing $C$.
  See also \cite{CY} for further results along these lines.  The inequality
(\ref{MT}) for $\gamma = 1$ and $\o =\o_{KE}$ corresponds to a reduction in the constant $\nu$
to $\nu=1-A_1$, with $C=e^{B_1}$.

\v

The exponent $\g$ found for (\ref{MT}) in \cite{T97}
was $\g=e^{-n}/(8n+8+e^{-n})$. In \cite{CTZ} it has been improved to $\g=1/(4n+5)$.
It was conjectured in \cite{T97} that the optimal
exponent should be $\g=1$. It has been shown in \cite{S04}
that this is indeed the case when the $\al$-invariant
\cite{T87}, \cite{TY} is strictly greater than ${n\over n+1}$.  Also, very recently,
it has been shown \cite{ZZ} that the conjecture holds
in the case of toric
manifolds for potentials
invariant under the compact torus action.

\v

In this paper, we prove Tian's conjecture in general.

\begin{theorem}\label{main} {\it Assume $c_1(X)>0$,  $X$ has no non-trivial holomorphic
vector fields and that $X$
admits a K\"ahler-Einstein metric $\o_{KE}$. Then
there exist positive constants $A$ and $B$ such that
\be\label{mt1} F_{\o_{KE}}(\phi) \ \geq \ AJ_{\o_{KE}}(\phi) \ - \ B \ \
\ee
for all  $\phi \in {\cal P} (X, \o_{KE})$.}
\end{theorem}

If $X$ does have non-trivial holomorphic vector fields, Theorem 1
can  be generalized as follows (cf. \cite{T97}, Theorem 4.4). Let
${\rm Aut}^0(X)$ be the component of the automorphism group of $X$
containing the identity, and $G={\rm Stab}(\o_{KE})\subset {\rm
Aut}^0(X)$ be the stabilizer of $\o_{KE}$.

\begin{theorem}
\label{main general}
Assume $c_1(X)>0$ and that $X$ admits a K\"ahler-Einstein metric
$\o_{KE}$. Assume that ${\rm Aut}^0(X)$ has a finite center. Then there
exist positive constants $A$ and $B$ such that
\be
\label{mtG}
F_{\o_{KE}}(\phi) \ \geq \ AJ_{\o_{KE}}(\phi) \ - \ B \ \
\ee
for all $G$-invariant potentials $\phi \in {\cal P} (X, \o_{KE})$. More
generally, if $G\subset {\rm Stab}(\o_{KE})$ is a closed subgroup
whose centralizer in ${\rm Stab}(\o_{KE})$ is finite, then the same
inequality (\ref{mtG}) holds for all $G$-invariant potentials $\phi$.
\end{theorem}

By the results of \cite{L}, \cite{Ma} and \cite{BM}, ${\rm Stab}(\o_{KE})$ is a maximal
compact subgroup of $\textrm{Aut}^0(X)$ and the complexification of its
Lie algebra  is the space of holomorphic vector fields.

\v

The proof of Theorem \ref{main} consists of two steps:

\bigskip
{\bf Step 1.}   Prove the inequality (\ref{MT}) for some $\gamma>0$ and for those $\phi$ satisfying
\bea \label{assume}
\textrm{osc}(\phi) \le K(J_{\o_{KE}}(\phi)+1),
\eea
for some constant $K$, where the $A_{\g}$ and $B_{\g}$ now depend on $K$.

\bigskip

{\bf Step 2.}  Remove the assumption (\ref{assume}) and prove the Moser-Trudinger inequality (\ref{MT}) with $\gamma=1$.

\bigskip

$\bullet$ Tian \cite{T97} has given a proof of the first step with $\gamma=e^{-n}/(8n+8+e^{-n})$.  His proof
makes use of the reverse continuity method of Bando-Mabuchi \cite{BM} and a modification of Bando's  smoothing by the K\"ahler-Ricci flow \cite{B}.  We have taken the opportunity to present a complete and somewhat more streamlined proof of Step 1 in sections 2-4.   There are some simplifications and improvements, and in particular we have eliminated the two Moser iteration arguments used in \cite{T97}.  As a consequence, $\gamma$ here can be taken to be arbitrarily close to $1/(2n-1)$ if $n>1$ and arbitrarily close to $\frac{1}{2}$ if $n=1$, an improvement over previously known results.

\bigskip

$\bullet$ It was shown by Tian-Zhu \cite{TZ} that the assumption (\ref{assume}) can be removed so that (\ref{MT}) holds with the $\g$ of Step 1.
 In section 5 we prove Step 2 by simultaneously removing the assumption and proving the full Moser-Trudinger inequality with $\g=1$, regardless of the value of $\g$ obtained in Step 1.
A reader who is familiar with the afore-mentioned result of \cite{T97} and wants to immediately see a proof of Theorem 1 may wish to skip to this section.

\bigskip

$\bullet$ Finally we show, also in section 5, how to prove Theorem 2 by making
 some minor
modifications to the proof of Theorem 1.

\bigskip
It may be interesting to have a version of Theorem 2 with the assumption of G-invariance replaced by a less restrictive orthogonality condition on $\phi$. If one considers the condition $\phi\perp\Lambda$, then one would have to show that the projection of $\psi_t$ on $\Lambda$ is controlled by the projection of $\psi_t$ on $\Lambda^\perp$, where $\psi_t$ are the potentials introduced in section 3 below. Note that the condition $\phi\perp\Lambda$ is different from the condition $e^{-\phi}\perp\Lambda$ introduced in \cite{A2}. It is known that (\ref{MT}) does not hold (for any exponent
$\gamma>0$) if the plurisubharmonicity condition $\omega_{KE}+{\sqrt{-1}\over 2}\ddb\phi>0$ is dropped, as there are counterexamples \cite{Sa} with $\phi\perp\Lambda$.
Other non-plurisubharmonic counterexamples with $\phi\perp\Lambda$ can also be deduced from the results of \cite{Be}.

\bigskip

{\bf Remark} \ The results above also immediately imply analogous inequalities
for the Mabuchi energy \cite{M} by the argument in \cite{T97} and the $E_1$ functional \cite{CT} by the argument in \cite{SW}.

\pagebreak[3]
\section{Proof of Step 1}
\setcounter{equation}{0}

We assume the existence of a K\"ahler-Einstein metric $\o_{KE}$,
and establish the inequality (\ref{MT}) for $F_{\o_{KE}}(\phi)$, for $\phi$ satisfying (\ref{assume}).
Fix a potential $\phi\in {\cal P}(X,\o_{KE})$, and
set $\o=\o_{KE}+{\sqrt{-1}\over 2}\ddb \phi$.
Consider the following complex Monge-Amp\`ere equation,
\bea
\label{MA1}
(\o+{\sqrt{-1}\over 2}\ddb\phi_t)^n
=e^{h_{\o}-t\phi_t}\o^n,
\qquad t\in[0,1].
\eea
By a theorem of Bando-Mabuchi \cite{BM}, since there are no non-trivial holomorphic
vector fields, a unique solution $\phi_t$ to this equation exists
for all $t\in [0,1]$, and $\o_{\phi_1}=\o_{KE}$.
In particular $\phi_1$ and $-\phi$ differ by a constant.
To exploit this equation, we require another functional $I_\o(\phi)$,
defined by
\bea
\label{I}
I_\o(\phi)
=
{\sqrt{-1}\over 2V}\sum_{j=0}^{n-1}\int_X
\pl\phi\wedge\bar\pl\phi\wedge \o^j\wedge \o_{\phi}^{n-1-j}
=
{1\over V}\int_X \phi(\o^n-\o_{\phi}^n).
\eea
From the explicit expressions for $I_\o$ and $J_\o$, it is easy to see that
$I_\o$, $J_\o$ are always positive, and that
\bea \label{ineqIJ}
{1\over n}J_\o \leq {1\over n+1}I_\o\leq J_\o.
\eea
The derivatives of these functionals $I_\o$, $J_\o$ along a general path $\phi_t$ are given by
\bea
{d \over dt} J_\o (\phi_t) & =& {1\over V}\int_X\dot\phi_t(\o^n-\o_{\phi_t}^n), \\
{d \over dt} I_\o (\phi_t) & =& {1\over V}\int_X\dot\phi_t(\o^n-\o_{\phi_t}^n)-
{1\over V}\int_X\phi_t {d \over dt} (\o_{\phi_t}^n).
\eea
As a first consequence of these formulas, we find
\bea \label{eqnIJ1}
{d \over dt} (I_\o-J_\o)(\phi_t)
=
-{d\over dt} \left( {1\over V}\int_X \phi_t\o_{\phi_t}^n \right)+{1\over V}\int_X\dot\phi_t\o_{\phi_t}^n.
\eea
When $\phi_t$ satisfies the Monge-Amp\`ere equation (\ref{MA1}),
we have $\int_X(t\dot\phi_t+\phi_t)\o_{\phi_t}^n=0$, which follows from differentiating the equation
$\int_X e^{h_{\o} - t \phi_t} \omega^n = V.$
Substituting for $\int_X\dot\phi_t\o_{\phi_t}^n$
in (\ref{eqnIJ1}) and integrating in $t$, we obtain \cite{Di} for all $t$
\bea \label{eqnDing}
-\frac{1}{V} \int_X\phi_t\o_{\phi_t}^n=(I_\o-J_\o)(\phi_t)
-{1\over t}\int_0^t(I_\o-J_\o)(\phi_s)ds,
\eea
or equivalently, using the expression (\ref{I}) for $I_\o(\phi_t)$,
\bea \label{eqnDing2}
F_\o(\phi_t)=-{1\over t}\int_0^t(I_\o-J_\o)(\phi_s)ds
-\log({1\over V}\int_Xe^{h_\o-\phi_t}\o^n).
\eea
Taking $t=1$, the last term on the right hand side is zero, and $F_\o(\phi_1)=-F_{\o_{KE}}(\phi)$,
so that
\bea
\label{FIJ}
F_{\o_{KE}}(\phi)=\int_0^1(I_\o-J_\o)(\phi_t)\ dt.
\eea
Next, the variational formulas for $I_\o$ and $J_\o$ also imply
\bea \label{eqnddtIminusJ}
{d \over dt} (I_\o-J_\o) (\phi_t)=-{1\over V}\int_X\phi_t {d \over dt} (\o_{\phi_t}^n)
=
-{1\over V}\int_X\phi_t(\Delta_t\dot\phi_t)\o_{\phi_t}^n,
\eea
where $\Delta_t$ is the scalar Laplacian with respect to the K\"ahler form $\o_{\phi_t}$.
Differentiating the Monge-Amp\`ere equation (\ref{MA1}) shows that
$\Delta_t\dot\phi_t=-t\dot\phi_t-\phi_t$, and hence
\bea \label{eqnIJ2}
{d \over dt} (I_\o-J_\o) (\phi_t)
=
\frac{1}{V} \int_X (-\Delta_t \dot\phi_t )(-\Delta_t-t)\dot\phi_t\,\o_{\phi_t}^n.
\eea
But the Monge-Amp\`ere equation (\ref{MA1}) also implies that
$\textrm{Ric}(\o_{\phi_t})=t\o_{\phi_t}+(1-t)\o \ge t\o_{\phi_t}$, and hence, by the Bochner-Kodaira formula,
$-\Delta_t-t$ is a positive operator. It follows that the right hand side of (\ref{eqnIJ2}) is nonnegative
and  $(I_\o-J_\o)(\phi_t)$ is a nondecreasing function of $t$.

\medskip

$\bullet$ We will also need the following consequence of the
fact that the functional $F_\o^0$
defined by
$F_\o^0(\phi)=J_\o(\phi)-{1\over V}\int_X\phi\o^n$ also
satisfies the cocycle property
$F_\o^0(\psi+\theta)=F_\o^0(\psi)+F_{\o_\psi}^0(\theta)$:
\bea
|J_\o(\phi_1)-J_\o(\phi_0)|
\leq {\rm osc}(\phi_1-\phi_0), \quad \textrm{for } \phi_0, \phi_1 \in {\cal
P}(X, \o).
\eea
Indeed,
\bea \nonumber
J_\o(\phi_1)-J_\o(\phi_0) &=&
 \frac{1}{V} \I_X(\phi_1-\phi_0)\o^n - \ F_{\o_{\phi_1}}^0(\phi_0-\phi_1) \\ \label{J}
& \leq & \frac{1}{V}  \I_X(\phi_1-\phi_0)\o^n + \frac{1}{V} \I_X
(\phi_0-\phi_1)\o_{\phi_{1}}^n \le \textrm{osc}(\phi_1-\phi_0),
\eea
where we made use in the inequality of the fact that $J \geq 0$.
Interchanging $\phi_0$ and $\phi_1$ establishes the
desired estimate. There is a similar inequality for the functional  $(I_{\o}-J_{\o})$ which will also be useful later:
\bea \label{IminusJ}
|(I_{\o} - J_{\o})(\phi_1) - (I_{\o}- J_{\o})(\phi_0) | & \le & n \, {\rm osc}(\phi_1 - \phi_0).
\eea
To see this, note that adding
\bea \nonumber
I_{\o}(\phi_0) - I_{\o}(\phi_1) & = & {1 \over V} \int_X (\phi_0 - \phi_1) \o^n + {1\over V} \int_X (\phi_1 - \phi_0) \o_{\phi_1}^n + {1 \over V} \int_X \phi_0 (\o_{\phi_1}^n - \o_{\phi_0}^n),
\eea
to (\ref{J}) gives
\bea \nonumber
(I_{\o} - J_{\o})(\phi_0) - (I_{\o} - J_{\o})(\phi_1)
 & \le & {1 \over V} \int_X \phi_0 (\o_{\phi_1}^n - \o_{\phi_0}^n) \\
& = &  {1 \over V} \int_X (\phi_1 - \phi_0) (\o_{\phi_0} - \o) \wedge \sum_{j=0}^{n-1} \o_{\phi_0}^j \wedge \o_{\phi_1}^{n-1-j},\qquad
\eea
after integrating by parts, and this is bounded above by $n \, \textrm{osc}(\phi_1 - \phi_0)$.  Interchanging $\phi_1$ and $\phi_0$ gives (\ref{IminusJ}).

\bigskip
We return now to the proof of the inequality (\ref{MT}).
Since $(I_{\omega}-J_{\omega})(\phi_t)$ is nondecreasing in $t$, the identity (\ref{FIJ}) implies that, for any $t$ in $[0,1]$,
\bea
F_{\o_{KE}}(\phi)\geq (1-t)(I-J)_{\o}(\phi_t)
\geq
{1\over n}(1-t)J_{\o}(\phi_t),
\eea
where we have made use of (\ref{ineqIJ}).
Since $|J_{\o}(\phi_t)-J_{\o}(\phi_1)|\leq {\rm osc}(\phi_t-\phi_1)$
and $J_{\o}(\phi_1)=J_{\o_{KE}}(\phi)$, we deduce that
\bea
\label{t}
F_{\o_{KE}}(\phi)
\geq {1\over n}(1-t)J_{\o_{KE}}(\phi)
-{1\over n}(1-t)\textrm{osc}(\phi_t - \phi_1).
\eea

\bigskip
$\bullet$ The main problem is then to estimate $\textrm{osc}(\phi_t - \phi_1)$.  It is of course sufficient to estimate $||\phi_1-\phi_t||_{C^0}$.
The original Monge-Amp\`ere equation (\ref{MA1}) is written
in terms of $\o_{\phi_t}\equiv\o+{\sqrt{-1}\over 2}\ddb \phi_t$,
so that $\o$ is the reference metric. Since the key estimates which we need
take place
near $t=1$, and since $\o$ depends on the given potential
$\phi\in {\cal P}(X,\o_{KE})$, it is preferable to rewrite (\ref{MA1})
as a Monge-Amp\`ere equation with $\o_{KE}$ as reference metric.  It is easily verified using (\ref{MA1}) that
\bea
\label{MA2}
\log {\o_{KE}^n
\over (\o_{KE}-{\sqrt{-1}\over 2}\ddb (\phi_1-\phi_t))^n}
+
(\phi_1-\phi_t)=(t-1)\phi_t,
\eea
For later use, we also
compare $(t-1)\phi_t$ with the potential function $h_{\o_{\phi_t}}$
for the Ricci curvature with the normalization (\ref{hnorm}).
We must have $h_{\o_{\phi_t}}=(t-1)\phi_t+c_t$
for some constant $c_t$.
The normalization condition (\ref{hnorm}) gives
\bea
\int_X\o_{\phi_t}^n
=
\int_Xe^{h_{\o_{\phi_t}}}\o_{\phi_t}^n
=
\int_X e^{(t-1)\phi_t+c_t}\o_{\phi_t}^n
\eea
from which it follows that $|c_t|\leq (1-t)||\phi_t||_{C^0}$.
In particular,
\bea
\label{h}
||h_{\o_{\phi_t}}||_{C^0}
\leq 2(1-t)||\phi_t||_{C^0}.
\eea

\bigskip
$\bullet$ The equation (\ref{MA2}) is of the form,
\bea
\label{MA3}
\log {\o_{KE}^n
\over (\o_{KE}-{\sqrt{-1}\over 2}\ddb \psi)^n}
+
\psi=h,
\eea
with $\psi = \phi_1 - \phi_t$ and $h= (t-1) \phi_t$.
As a first attempt to get estimates for $||\phi_1-\phi_t||_{C^0}$ in terms
of $1-t$ and $||\phi_t||_{C^0}$, we can apply the implicit function theorem.
The linearization of the left hand side of (\ref{MA3}) at $\psi=0$ is
$\delta\psi \mapsto \Delta_{KE}\delta\psi+\delta\psi$,
which is an invertible operator from $L_{k+2}^p(X)\to L_{k}^p(X)$
and from $C^{2,{\kappa}}(X)\to C^{\kappa}(X)$,
for any $p>1$, $k \ge 0$ and $0<\kappa<1$, when $X$ has no non-trivial holomorphic vector
fields. Thus there exist constants $C(\o_{KE})$, $\epsilon(\o_{KE})>0$
(depending on $p$, $k$, and $\kappa$)
so that if $\psi$ and $h$ satisfy the equation (\ref{MA3}), then
\bea
\label{IFT}
||h||_{C^{0,{\kappa}}} \leq \epsilon(\o_{KE})
\ &\Rightarrow& \ ||\psi||_{C^{2,{\kappa}}}\leq C(\o_{KE})\ ||h||_{C^{0,\kappa}}
\eea
and
\bea
||h||_{L_{k}^p} \leq \epsilon(\o_{KE})
\ &\Rightarrow& \ ||\psi||_{L_{k+2}^p}\leq C(\o_{KE})\ ||h||_{L_k^p}.
\eea
The above norms, and indeed any other unadorned norms appearing in this paper,
are understood to be taken with respect to the metric $\o_{KE}$.
Now for $p> 2n$ and $k=0$, we obtain, in view of the Sobolev imbedding theorem,
\bea
||h||_{C^0}\leq \epsilon(\o_{KE})
\ \Rightarrow\ ||\psi||_{C^0} \leq C(\o_{KE})\ ||h||_{C^0}.
\eea
 Applied to the equation (\ref{MA2}), we find
\bea
\label{imply}
(1-t)||\phi_t||_{C^0}
\leq\epsilon(\o_{KE})
\ \Rightarrow\ ||\phi_1-\phi_t||_{C^0}\leq C (\o_{KE})\ (1-t)||\phi_t||_{C^0}.
\eea
If indeed $(1-t)||\phi_t||_{C^0}\leq\epsilon(\o_{KE})$
for all $t$, then we pick say $t={1\over 2}$
in the inequality (\ref{t}), and the desired inequality
(\ref{mt1}) follows at once.
Otherwise, we let $t_0$ be the largest value possible
so that $(1-t)||\phi_t||_{C^0}<\e(\o_{KE})$ for all $t\in(t_0,1]$.
Clearly $t_0<1$, although it does depend on the potential
$\phi\in {\cal P}(X,\o_{KE})$. Assuming without loss of generality
that both $\epsilon(\o_{KE})$ and $1-t_0$ are small numbers, we can
deduce that ${1\over 2}||\phi_{t_0}||_{C^0}
\leq ||\phi_1||_{C^0}\leq 2||\phi_{t_0}||_{C^0}$,
and hence
\bea
F_{\o_{KE}}(\phi)
\ \geq\ {\epsilon(\o_{KE})\over 2n}{J_{\o_{KE}}(\phi)\over||\phi_1||_{C^0}}
-C'.
\eea
However, in general, $||\phi_1||_{C^0}$ is of the same size as $J_{\o_{KE}}$, so
this estimate does not improve over the already known fact that $F_{\o_{KE}}(\phi)$
is bounded from below.

\bigskip
$\bullet$ The idea is to improve on the range of $t$
on the left side of the implication in (\ref{imply}) so that, roughly speaking, $(1-t) \sim 1/\| \phi_{t} \|_{C^0}^{1-\g}$ for some $\g>0$, assuming $\| \phi_t \|_{C^0}$ is large.   More precisely, we have
the following lemma, which provides the key estimate for
$||\phi_1-\phi_t||_{C^0}$:

\bigskip
{\bf Lemma 1} {\it Let $p>2n$, $0<\kappa<1$ and $\alpha={p+\kappa-2\over p-1}$.
Then there exists a constant
$D=D(\o_{KE})>0$ depending on $p$, $\kappa$ so that
\bea
\label{phi1phit}
||\phi_1-\phi_t||_{C^0}
\leq 100(1-t)||\phi_t||_{C^0}
+
1,
\eea
for all $t\in [t_0,1]$, where $t_0$ (depending on $\phi$) is defined by}
\bea
\label{t0}
(1-t_0)^{1-\alpha}(1+2(1-t_0)||\phi_{t_0}||_{C^0})^\alpha
  =
\sup_{t\in [t_0,1]}
(1-t)^{1-\alpha}(1+2(1-t)||\phi_t||_{C^0})^\alpha
=
D. \qquad
\eea

\bigskip
Given Lemma 1, we can complete Step 1. Assume that $p > 3 - 2\kappa$,
so that $\alpha > 1/2$.
 Notice that this holds automatically if $n>1$.
In the inequality (\ref{t}),
we set $t=t_0$, where $t_0$ is chosen as in Lemma 1.
The definition (\ref{t0}) of $t_0$ implies that
\bea
2^\al(1-t_0)||\phi_{t_0}||_{C^0}^{\al}
\leq D(\o_{KE}).
\eea
Then since $1/\alpha < 2$,
\bea
(1-t_0)^2||\phi_{t_0}||_{C^0}\leq D(\o_{KE})^{1\over \al}.
\eea
Hence,  from the inequality (\ref{phi1phit}),
\bea
(1-t_0)||\phi_1-\phi_{t_0}||_{C^0}
\leq
100(1-t_0)^2||\phi_{t_0}||_{C^0}+1
\leq C_1,
\eea
where $C_1$ is another constant depending only on $\o_{KE}$, $p$ and $\kappa$.
This handles the second term on the right hand side of (\ref{t}).
To handle the first term on the right hand side, we begin by noting
that we may assume that $1-t_0<\delta$ for some fixed small
constant $\delta>0$, otherwise the stronger
inequality $F_{\o_{KE}}(\phi)\geq {\delta\over n}J_{\o_{KE}}(\phi)-{2\over n}C_1$
already holds. In particular, we may assume that
$(1-t_0)||\phi_{t_0}||_{C^0}\geq {1\over 2}$, for otherwise
the definition (\ref{t0}) implies that
$2^{-{\al\over{1-\alpha}}}D(\o_{KE})^{1\over {1-\alpha}}< 1-t_0$.
When $(1-t_0)||\phi_{t_0}||_{C^0}\geq {1\over 2}$, the definition (\ref{t0})
implies
\bea
D(\o_{KE})
\leq 2^{2\al} (1-t_0)||\phi_{t_0}||_{C^0}^{\al},
\eea
giving
\bea
1-t_0\geq {C_2\over ||\phi_{t_0}||_{C^0}^{\al}}.
\eea
The inequality (\ref{t}) implies then
\bea
F_{\o_{KE}}(\phi)
\geq C_2{J_{\o_{KE}}(\phi)
\over
||\phi_{t_0}||_{C^0}^{\al}}-{2\over n}C_1.
\eea
However, for $\delta\leq {1\over 200}$, the inequality (\ref{phi1phit})
implies that
\bea
||\phi_{t_0}||_{C^0}\leq 2||\phi_1||_{C^0}+2,
\eea
so we obtain
\bea
\label{2}
F_{\o_{KE}}(\phi)
\geq C_2{J_{\o_{KE}}(\phi)
\over
(2 ||\phi_{1}||_{C^0}+2)^{\al}}-{2\over n}C_1.
\eea
Since
\bea
\int_X e^{t\phi_t}\o_{\phi_t}^n=\int_X e^{h_{\o}}\o^n=\int_X\o^n
=
\int_X\o_{\phi_t}^n,
\eea
$\phi_t$ must change signs.  Hence  $||\phi_t||_{C^0}\leq  {\rm osc}(\phi_t)$.
 Furthermore ${\rm osc}\,\phi_1={\rm osc}\,\phi$
since $\phi_1$ and $-\phi$ differ by a constant. Thus we have shown that
\bea
\label{MTK}
F_{\o_{KE}} \geq C_2 \ {J_{\o_{KE}}(\phi)
\over (2\, {\rm osc}\,\phi+2)^{\al}}-{2\over n}C_1.
\eea
The desired inequality (\ref{MT}) with $\g=1-\alpha
={1-\kappa\over p-1}$
and all $\phi\in{\cal P}(X,\o_{KE})$
satisfying the additional condition (\ref{assume}).
If $n>1$,  $p$ and $\kappa$ can be taken arbitrarily close to $2n$
and to $0$ respectively and it follows
that $\g$ can be taken arbitrarily close to ${1\over 2n-1}$.  If $n=1$ then
we require $p>3-2\kappa$ giving a $\g$ arbitrarily close
to $\frac{1}{2}$.  This completes Step 1, assuming Lemma 1.

\section{Proof of Lemma 1}
\setcounter{equation}{0}

The idea in the proof of Lemma 1 is to smooth out $\phi_1-\phi_t$ by
the K\"ahler-Ricci flow. More precisely, for each $t$, consider the
normalized K\"ahler-Ricci flow $u_{s,t}$ in time $s$ with initial data $u_{0,t} = 0$,
\bea
\label{KRflow} {\pl u_{s,t}\over\pl s} = \log
{(\o_{\phi_t}+{\sqrt{-1}\over 2}\ddb u_{s,t})^n\over\o_{\phi_t}^n}
+u_{s,t}-h_{\o_{\phi_t}}.
\eea
Write $u_t$ for $u_{1,t}$,
and consider the K\"ahler form
\bea
\o_{\phi_t+u_t}=
\o+{\sqrt{-1}\over 2}\ddb (\phi_t+u_t)
=
\o_{KE}-{\sqrt{-1}\over 2}\ddb(\phi_1-\phi_t-u_t).
\eea
There exists then a constant $a_t$ so that
\bea
\log \left( {\o_{KE}^n\over (\o_{KE}-{\sqrt{-1}\over 2}\ddb(\phi_1-\phi_t-u_t))^n}\right)
+\phi_1-\phi_t-u_t-a_t=h_{\o_{\phi_t+u_t}},
\eea
which can be seen easily by applying $\sqrt{-1} \ddb$ to both sides.
Set now
\bea
\psi_t=\phi_1-\phi_t-u_t-a_t.
\eea
Then the equation (\ref{MA3}) holds with
$\psi=\psi_t$, $h=h_{\o_{\phi_t+u_t}}$. The following ``smoothing lemma"
shows how the bound on $h_{\o_{\phi_t+u_t}}$ improves on that of $h_{\o_{\phi_t}}$:

\bigskip
{\bf Lemma 2} {\it Let $u_t = u_{1,t}$ be defined as in (\ref{KRflow}).
\begin{enumerate}
\item[(a)]
We have for all $t\in[0,1],$
\bea
||u_t||_{C^0}\leq 3||h_{\o_{\phi_t}}||_{C^0}.
\eea
\item[(b)] Let $A>0$
and $t_1 \in [0,1]$. Assume that
\bea
\label{A}
A^{-1}\,\o_{KE}\leq \o_{\phi_t+u_t} \leq A\, \o_{KE},
\eea
for all $t \in [t_1,1]$.  Then for any $p>2n$ and $0<\kappa<1$,
there exists a constant $B_A$ depending only on $p$, $A$, and $\o_{KE}$ so
that for $t\in [t_1,1]$,
\bea
\label{hsmoothed}
||h_{\o_{\phi_t+u_t}}||_{C^{0,{\kappa}}}
\leq B_A(1+||h_{\o_{\phi_t}}||_{C^0})^{\al}(1-t)^{1-\alpha}.
\eea
where
$\alpha={p+\kappa-2\over p-1}$.  Here, as before, the norms are with respect
to $\o_{KE}$.
\end{enumerate}}

\bigskip
This lemma follows from estimates for the K\"ahler-Ricci flow and will be proved
in the next section. Assuming it for the moment, we complete now the proof of Lemma 1.
The next step is to establish the following lemma:

\bigskip
{\bf Lemma 3} {\it Fix $p>2n$, $0<\kappa<1$,
and choose $D(\o_{KE})={\e(\o_{KE})\over 4(B_2+1)(C+1)(\e+1)}$,
where $\e=\e(\o_{KE})$, $C=C(\o_{KE})$ are chosen as in (\ref{IFT}), and $B_2$ is
defined as in Lemma 2
with $A=2$.
Choose $t_0$ as in (\ref{t0}). Then for $t\in [t_0,1]$, we have,
}
\bea
||\psi_t ||_{C^{2,{\kappa}}} <  {1\over 4}.
\eea

\bigskip
{\it Proof of Lemma 3.} We assume the contrary. Since $\psi_1=0$,
there exists $t_1\in [t_0,1)$ so that
\bea
||\psi_{t_1}||_{C^{2,{\kappa}}}={1\over 4}
\qquad
{\rm and}
\qquad ||\psi_t||_{C^{2,{\kappa}}}<{1\over 4}
\quad {\rm if}
\quad t>t_1.
\eea
Since the operator norm of any symmetric matrix
is always smaller than its Hilbert-Schmidt norm,
we have for $t\in [t_1,1]$,
\bea
-{1\over 4}\o_{KE}
\leq {\sqrt{-1} \over 2} \ddb\psi_t\leq {1\over 4}\o_{KE},
\eea
and thus
\bea
\label{obound}
{1\over 2}\o_{KE}
\leq \o_{KE}-{\sqrt{-1}\over 2}\ddb\psi_t
=
\o_{\phi_t+u_t}
\leq 2\,\o_{KE}
\qquad {\rm for}\ t\in[t_1,1].
\eea
Thus the inequality (\ref{h}) and Lemma 2 imply that for $t\in [t_1,1]$,
\bea
||h_{\o_{\phi_t+u_t}}||_{C^{0,{\kappa}}}
&\leq& B_2(1+||h_{\o_{\phi_t}}||_{C^0})^{\al}(1-t)^{1-\alpha}
\nonumber\\
&\leq& B_2(1+2(1-t)||\phi_t||_{C^0})^\al(1-t)^{1-\alpha}
\nonumber\\
&\leq&
{\e B_2\over 4(B_2+1)(C+1)(\e+1)}.
\eea
This bound is strictly less than $\e$,
and thus the bound (\ref{IFT}) from the implicit function theorem implies for $t\in [t_1,1]$
\bea
||\psi_t||_{C^{2,{\kappa}}}
\leq C\ ||h_{\o_{\phi_t+u_t}}||_{C^{0,{\kappa}}}
\le  C\,{\e B_2\over 4(B_2+1)(C+1)(\e+1)}
\,<{1\over 4}.
\eea
This contradicts $||\psi_{t_1}||_{C^{2,{\kappa}}}={1\over 4}$,
and completes the proof of the lemma.  Q.E.D.

\bigskip
{\it Proof of Lemma 1.} Still assuming Lemma 2, we can give now the proof of
Lemma 1. Since $\phi_1-\phi_t=\psi_t+u_t+a_t$, and since the bound $\| \psi_t \|_{C^0} < 1$
is provided by Lemma 3, it will suffice to estimate $||u_t||_{C^0}$
and $|a_t|$. Using the estimate (\ref{h}) for $h_{\o_{\phi_t}}$,
Lemma 2 implies at once $||u_t||_{C^0}\leq 6(1-t)||\phi_t||_{C^0}$.
Next, the Monge-Amp\`ere equation for $\psi_t$ and the normalization
of $h_{\o_{\phi_t+u_t}}$ give
\bea
\int_X e^{\psi_t} \o_{KE}^n
=
\int_X e^{h_{\o_{\phi_t+u_t}}}\o_{\phi_t+u_t}^n=V.
\eea
On the other hand, rewriting the Monge-Amp\`ere equation (\ref{MA2}) in the form
$\o_{KE}^n=e^{t\phi_t-\phi_1}\o_{\phi_t}^n$, we also have
\bea
\int_X e^{\psi_t}\o_{KE}^n
=
\int_X e^{\psi_t+t\phi_t-\phi_1}\o_{\phi_t}^n
=
\int_Xe^{-a_t-u_t+(t-1)\phi_t}\o_{\phi_t}^n.
\eea
Comparing the two values, we obtain
\bea
|a_t|\leq ||u_t||_{C^0}+(1-t)||\phi_t||_{C^0}
\leq 7(1-t)||\phi_t||_{C^0}.
\eea
We can return now to the estimate of $||\phi_1-\phi_t||_{C^0}$,
\bea
||\phi_1-\phi_t||_{C^0}
\leq |a_t|+||u_t||_{C^0}+||\psi_t||_{C^0},
\eea
and the desired estimate for $||\phi_1-\phi_t||_{C^0}$ follows immediately.
The proof of Lemma 1 is complete.  Q.E.D.

\section{Smoothing by the K\"ahler-Ricci Flow}
\setcounter{equation}{0}

We give now the proof of Lemma 2. It is here that the improvement
of the exponent in Step 1 over previous results takes place. For convenience,
introduce the following notation for the K\"ahler-Ricci flow,
\bea \label{eqnKR}
{\pl u\over\pl s}
=
\log {(\eta_0+{\sqrt{-1}\over 2}\ddb u)^n\over\eta_0^n}
+
u-h_{\eta_0},
\qquad u(0)=0,
\eea
$\eta_s=\eta_0+{\sqrt{-1}\over 2}\ddb u$, and $h_s =
h_{\eta_s}$. Then $h_{s}=-\dot u+c_s=\tilde h_s+c_s$
for some constant $c_s$ with $c_0=0$.  We will use a subscript $s$  to indicate  objects  that are defined
with respect to the metric $\eta_s$ (e.g. $\Delta_s$,
$\nabla_s$, $| \cdot |_s$).
 The following lemma is well-known
\cite{B}, but we include its proof for completeness:

\bigskip
{\bf Lemma 4}. {\it We have the following inequalities:
}
\bea
&{\rm (a)}&\qquad\qquad\qquad ||\dot u||_{C^0}\leq e^s||h_0||_{C^0}
\nonumber\\
&{\rm (b)}&\qquad\qquad\qquad \sup_X (|h_s|^2+s|\nabla h_s|_s^2) \leq 4 e^{2s}||h_0||_{C^0}^2
\nonumber\\
&{\rm (c)}&\qquad\qquad\qquad  e^{-s}\Delta_sh_s\geq \Delta_0h_0.
\eea

\bigskip
{\it Proof of Lemma 4.} Differentiating the K\"ahler-Ricci flow equation gives
\bea
{\pl\over\pl s}\dot u=\Delta_s\dot u+\dot u,
\eea
and hence $||\dot u||_{C^0}\leq e^s||h_0||_{C^0}$, giving (a).
Similarly, the flows for $\tilde h_s^2$ and $|\nabla \tilde h_s|_s^2$ are,
\bea
{\pl\over\pl s}\tilde h_s^2&=&
\Delta_s \tilde h_s^2-2|\nabla \tilde h_s|_s^2+2\tilde h_s^2
\nonumber\\
{\pl\over\pl s}|\nabla\tilde  h_s|_s^2
&=&
\Delta_s|\nabla\tilde  h_s|^2-|\nabla \bar\nabla\tilde  h_s|^2_s
-
|\nabla_s \nabla\tilde  h_s|^2_s+|\nabla\tilde  h_s|_s^2.
\eea

and hence,
\bea
{\pl\over\pl s}(\tilde h_s^2+s|\nabla\tilde  h_s|_s^2)
\leq
\Delta_s (\tilde h_s^2+s|\nabla\tilde  h_s|_s^2)
+2(\tilde h_s^2+s|\nabla\tilde  h_s|_s^2).
\eea
The maximum principle then implies
\bea\label{tilde}
\sup_X\big(e^{-2s}(\tilde h_s^2+s|\nabla\tilde  h_s|_s^2)\big)\leq ||h_0||_{C^0}^2 .
\eea

From the normalization condition for $h_s$ and (\ref{eqnKR}) we see that
\bea \label{nheta}
\int_X e^{h_{\eta_0} + c_s -u} \eta_0^n = V\ .
\eea
Part (a) implies $|c_s| \le e^s \| h_0 \|_{C^0}$ which, when combined with
(\ref{tilde}), yields (b).

\v

Finally, the flow for $\Delta_s \dot u$ is given by
\bea
{\pl\over\pl s}(\Delta_s \dot u)
=
\Delta_s^2\dot u+\Delta_s \dot u
-
|\nabla\bar\nabla\dot u|_s^2
\leq
\Delta_s^2\dot u+\Delta_s \dot u,
\eea
and (c) follows also from the maximum principle. Q.E.D.

\bigskip

{\bf Lemma 5} {\it Assume now that $\eta_0=\o_{\phi_t}$,
$h_0=h_{\o_{\phi_t}}$. Let $v=h_1-{1\over V}\int_X h_1\,\eta_1^n$
and assume that
\bea
\label{metricequiv2}
A^{-1}\o_{KE} \leq \eta_1 \leq A\,\o_{KE}.
\eea
Then for any $p>2n$,
there exists $C>0$, depending only on $\o_{KE}$, $A$ and $p$ so that}
\bea
||v||_{C^0}
\leq C\,||h_0||_{C^0}^{p-2\over p-1}(1-t)^{1\over p-1}.
\eea

\bigskip
{\it Proof of Lemma 5}. Lemma 4 shows that
\bea
||v||_{C^0}\leq  4 e ||h_0||_{C^0}.
\eea
On the other hand, we also have
\bea
\int_X|\nabla v|_1^2\eta_1^n
&=&-\int_X v (\Delta_1v) \cdot \eta_1^n
=
\int_X(v-\inf_X v)(-\Delta_1v)\cdot\eta_1^n
\nonumber\\
&\leq&
\int_X (v- \inf_X v)\sup_X(-\Delta_1v)\cdot\eta_1^n
\nonumber\\
&\leq&
2V ||v||_{C^0}\,\sup_X(-\Delta_1v).
\eea
Recall that $h_0=h_{\o_{\phi_t}}$ and thus
$\textrm{Ric}(\eta_0)>t\eta_0$. This implies that $\Delta_0h_0\geq -n(1-t)$,
and hence, in view of Lemma 4,
\bea
-\Delta_1h_1\leq ne(1-t).
\eea
Substituting in the previous inequality gives
\bea
\label{gradient}
\int_X|\nabla v|_1^2\eta_1^n
\leq 2V ne||v||_{C^0}(1-t).
\eea
Let $p>2n$. Then for some constants $C_i$ depending only on $\o_{KE}$,
$A$ and $p$,
\bea
||v||_{C^0}^p
&\leq& C_0 \,\left(\int_X |v|^p\eta_1^n
+
\int_X|\nabla v|_1^p\,\eta_1^n\right)
\nonumber\\
&\leq&
C_0\,
\left(||v||_{C^0}^{p-2}\int_X v^2\eta_1^n
+
(2e ||h_0||_{C^0})^{p-2}\int_X|\nabla v|_1^2\eta_1^n\right)
\nonumber\\
&\leq&
C_1\,||h_0||_{C^0}^{p-2}\int_X |\nabla v|_1^2\eta_1^n,
\eea
where we have used the Sobolev inequality, the Poincar\'e inequality and applied (b) of Lemma 4.
Note that the constants in the Sobolev and Poincar\'e inequalities depend only on
$\o_{KE}$, since the metric $\eta_1$ is equivalent
to $\o_{KE}$ in view of (\ref{metricequiv2}).
Together with the inequality (\ref{gradient}), this gives
\bea
||v||_{C^0}^p
\leq C_2(1-t)||h_0||_{C^0}^{p-2}||v||_{C^0},
\eea
which is the inequality to be proved. Q.E.D.

\bigskip
{\it Proof of Lemma 2.} We can give now the proof of Lemma 2.
Recall that $\eta_0=\o_{\phi_t}$ and $h_0=h_{\o_{\phi_t}}$.
From (a) of Lemma 4, it follows that
$-e^s||h_0||_{C^0}\leq \dot u\leq e^s||h_0||_{C^0}$,
and integrating from $0$ to $1$, we obtain the first
statement in Lemma 2. Next, fix $0<\kappa<1$
and let $x,y\in X$.  Let $d_1(x,y)$ be the distance between $x$ and $y$ with
respect to the metric $\eta_1$.
If $d_1(x,y)\geq (1+||h_0||_{C^0})^{-\mu} (1-t)^{\delta}$,
for some fixed $0<\delta,\mu\leq 1$ to be chosen later,
then Lemma 5 implies
\bea
{|h_1(x)-h_1(y)|\over
d_1(x,y)^\kappa}
&\leq& 2||v||_{C^0}(1+||h_0||_{C^0})^{{\mu\kappa}}(1-t)^{-\delta\kappa}
\nonumber\\
&\leq&
C (1+||h_0||_{C^0})^{{p-2\over p-1}+{\mu\kappa}}(1-t)^{{1\over p-1}-{\delta\kappa}}.
\eea
If $d_1(x,y) < (1+||h_0||_{C^0})^{-\mu}(1-t)^{\delta}$ instead,
then we use the estimate $||\nabla h_1||_{C^0(\eta_1)}
\leq 2e||h_0||_{C^0}$ from (b) of Lemma 4 to write
\bea
{|h_1(x)-h_1(y)|\over
{d_1(x,y)^\kappa}}
&\leq &
d_1(x,y)^{1-\kappa}||\nabla h_1||_{C^0(\eta_1)}
\nonumber\\
&\leq&
C\,(1+||h_0||_{C^0})^{-{\mu(1-\kappa)}+1}(1-t)^{\delta(1-\kappa)}.
\eea
The optimal choice of $\delta,\mu$ is when the exponents of $(1-t)$
and $||h_0||_{C^0}$ match in the two inequalities,
which is $\delta=\mu={1\over p-1}$. We get then in all
cases,
\bea
\label{Ckappa}
{|h_1(x)-h_1(y)|\over
{d_1(x,y)^\kappa}}
\leq
C\,(1+||h_0||_{C^0})^{\alpha} (1-t)^{1- \alpha}
\eea
with $\alpha={p+\kappa-2\over p-1}$.
Finally, $\int_Xe^{h_1}\eta_1^n=\int_X\eta_1^n$ implies
that $h_1$ changes signs, and thus
\bea
||h_1||_{C^0}
\leq
{\rm osc}\,(h_1)
=
{\rm osc}\,(v)
\leq
2||v||_{C^0}
\leq
C(1-t)^{1\over p-1}||h_0||_{C^0}^{p-2\over p-1}.
\eea
This is a better bound than the one in (\ref{Ckappa}),
so we obtain the
desired bound for $||h_1||_{C^{0,\kappa}}$. The
proof of Lemma 2, and hence of Step 1, is complete. Q.E.D.

\section{Proofs of the Main Theorems}
\setcounter{equation}{0}

\bigskip

In this section we complete Step 2 to prove Theorem 1 and then show how the
proof of Theorem 1 can be modified to give a proof of Theorem 2.

\bigskip
\noindent
{\it Proof of Theorem 1.} \,
We need to prove that for all $\phi \in {\mathcal P} (X, \o_{KE})$,
\bea \label{tc}
F_{\o_{KE}} (\phi) \ge A J_{\o_{KE}}(\phi) - B,
\eea
for some $A,B>0$.  As before, let $\omega = \o_{KE} + {\sqrt{-1} \over 2}
\ddb \phi$ and let $\phi_t$ be the solution of (\ref{MA1}).
 We first need the observation of Tian-Zhu \cite{TZ} that the potential $\phi_t - \phi_1$ satisfies the condition (\ref{assume}),
\bea
\label{eqnTZ1}
\textrm{osc}(\phi_t - \phi_1) & \le &  K (J_{\o_{KE}} (\phi_t - \phi_1)+1),
\eea
for $t \in [1/2,1]$ with a constant $K$ depending only on $\o_{KE}$.  We will then prove (\ref{tc}) by showing that there exist positive constants $C$ and $c$ depending only on $\o_{KE}$ such that for some $t$ with $1-t \ge c$,
\bea \label{Jbound}
J_{\o_{KE}}(\phi_{t} - \phi_1) \le C.
\eea
Recall from section 2 that we have the inequality (\ref{t}),
\bea
F_{\o_{KE}}(\phi) \ge \frac{1}{n} (1-t) J_{\o_{KE}}(\phi) - \frac{1}{n} ( 1-t)
\textrm{osc}(\phi_t - \phi_1),
\eea
for any $t \in [0,1]$.  Then this together with
 (\ref{eqnTZ1}) and (\ref{Jbound}) will give (\ref{tc}).

\bigskip

To prove (\ref{eqnTZ1}) observe that $\textrm{Ric} (\o_{\phi_t}) \ge \frac{1}{2} \o_{\phi_t}$ when $t $ is in the interval $[1/2,1]$.  It follows (see \cite{Si}, for example) that the Green's function for $\o_{\phi_t}$ is uniformly bounded from below.  Since $\o_{KE} + {\sqrt{-1} \over 2} \ddb (\phi_t - \phi_1) = \o_{\phi_t}>0$ and $\o_{\phi_t} + {\sqrt{-1} \over 2} \ddb(\phi_1 - \phi_t) = \o_{KE}>0$ we have
\bea
 - \Delta_{\o_{KE}} (\phi_t - \phi_1) <n \quad \textrm{and} \quad - \Delta_t (\phi_t - \phi_1) > - n.
\eea
Then using the Green's functions for $\o_{KE}$ and $\o_{\phi_t}$ there exists  $C_0=C_0(\o_{KE})$ such that for all $x,y$ in $X$,
\bea
(\phi_t- \phi_1)(x) & \le & \frac{1}{V} \int_X (\phi_t - \phi_1) \o_{KE}^n + C_0 \\
(\phi_t - \phi_1)(y) & \ge & \frac{1}{V} \int_X (\phi_t - \phi_1) \o_{\phi_t}^n - C_0.
\eea
Hence, using (\ref{I}) and (\ref{ineqIJ}),
\bea \nonumber
\textrm{osc}(\phi_t - \phi_1) & \le & \frac{1}{V} \int_X (\phi_t - \phi_1) (\o_{KE}^n - \o_{\phi_t}^n) + 2C_0 \\
& = &  I_{\o_{KE}} (\phi_t - \phi_1) + 2C_0 \le K (J_{\o_{KE}}(\phi_t - \phi_1)+1),
\eea
as required.

\bigskip

We will now give a proof of (\ref{Jbound}). Making use of the cocycle condition for $F$, the equations (\ref{MA1}), (\ref{eqnDing}), (\ref{eqnDing2}) and the concavity of the log function, we have for $t \in (0,1]$,
\bea \nonumber
F_{\o_{KE}} (\phi_t - \phi_1) & = & F_{\o}(\phi_t) - F_{\o}(\phi_1) \\ \nonumber
& = &  - {1 \over t} \int_0^t (I_{\o}- J_{\o})(\phi_s)ds + \int_0^1 (I_{\o} - J_{\o})(\phi_s) ds - \log\left( {1 \over V} \int_X e^{(t-1)\phi_t} \o_{\phi_t}^n \right) \\ \nonumber
& \le & - {1 \over t} \int_0^t (I_{\o}- J_{\o})(\phi_s)ds + \int_0^1 (I_{\o} - J_{\o})(\phi_s) ds + {(1-t) \over V} \int_X \phi_t \, \o_{\phi_t}^n \\ \nonumber
& = & - {1 \over t} \int_0^t (I_{\o} - J_{\o})(\phi_s)ds + \int_0^1 (I_{\o} - J_{\o}) (\phi_s)ds - (1-t) (I_{\o} - J_{\o})(\phi_t) \\ \nonumber
&& \mbox{} + {(1-t) \over t} \int_0^t (I_{\o} - J_{\o})(\phi_s) ds \\ \nonumber
& = & \int_t^1 (I_{\o} - J_{\o} )(\phi_s)ds - (1-t) (I_{\o} - J_{\o})(\phi_t) \\  \nonumber
& \le & (1-t) \left( (I_{\o} - J_{\o})(\phi_1) - ( I_{\o} - J_{\o})(\phi_t) \right) \\
& \le & n(1-t) \, \textrm{osc}(\phi_t - \phi_1),
\eea
where for the last line we have used  inequality (\ref{IminusJ}) and in the
previous line, we have used the fact that $(I_{\o}-J_{\o})(\phi_s)$ is nondecreasing
in $s$.  From  (\ref{eqnTZ1}) we see that for $t \in [1/2,1]$,
\bea
F_{\o_{KE}} (\phi_t - \phi_1) \le  n (1-t) K (J_{\o_{KE}}(\phi_t-\phi_1) +1).
\eea
But from Step 1 we have
\bea
 A_{\g} J_{\o_{KE}}(\phi_t - \phi_1)^{\g} - B_{\g} \le F_{\o_{KE}} (\phi_t - \phi_1),
\eea
and combining this with the preceding inequality gives
\bea \label{ineqJ}
J_{\o_{KE}}(\phi_t - \phi_1)^{\g} \left( A_{\g} - nK(1-t) J_{\o_{KE}}(\phi_t - \phi_1)^{1-\g} \right) \le n(1-t)K + B_{\g}.
\eea
We can suppose that there exists a $t' \in [1/2,1]$ with
\bea \label{eqnt'}
nK(1-t') J_{\o_{KE}}(\phi_{t'} - \phi_1)^{1-\g} = A_{\g}/2.
\eea
For, if not then we must have
\bea
nK(1-t) J_{\o_{KE}}(\phi_t - \phi_1)^{1-\g} < A_{\g}/2,
\eea
for all $t \in [1/2,1]$ and hence for $t=1/2$.  It would follow that $J_{\o_{KE}}(\phi_{1/2} - \phi_1) \le C$, and we would be done.  Otherwise, from (\ref{ineqJ}) we see that
$J_{\o_{KE}}(\phi_{t'} - \phi_1) \le C$ and from (\ref{eqnt'}) that $1-t' \ge c >0$ and we are finished. Q.E.D.
\bigskip

\noindent
{\it Proof of Theorem 2.} \, The assumption of no non-trivial holomorphic vector fields
is used
in the proof of Theorem 1 for two applications of the implicit function
theorem.  The first application is to obtain a solution $\phi_t$ of (\ref{MA1})
and the second is to obtain the estimate of Lemma 3.  In both cases, the
key fact is that the operator $f \mapsto \Delta_{\o_{KE}} f + f$ is an
invertible operator from $C^{2, \kappa}(X)$ to $C^{\kappa}(X)$. The
kernel $\Lambda$ of this operator  will be nontrivial if there are
holomorphic vector fields. Let $G$ be as in Theorem 2, and let $\lambda_1,\cdots,\lambda_M$
be a basis of $\Lambda$, and define for each $g\in G$ the matrix $a_{ij}(g)$
by $\rho(g)\lambda_i=a_{ij}(g)\lambda_j$, where $\rho$ is the action of $G$
on $\Lambda$. Then if $f$ is $G$-invariant, we have
\bea
V_i\equiv\int_X f \lambda_i \, \o_{KE}^n = a_{ij}(g)\int_X f\lambda_j\,\o_{KE}^n
=a_{ij}(g)V_j,
\eea
which means that the vector $V_i$ defined by the left hand side is fixed under $G$.
Since $G$ has by assumption no non-zero fixed vector, the vector $V_i$ must be $0$,
which shows that $f$ is orthogonal to $\Lambda$.
Hence replacing $C^{2, \kappa}(X)$ and $C^{\kappa}(X)$ with their $G$-invariant
counterparts, we see that the operator is still invertible in this case as long
as the objects involved are $G$-invariant.  But
  $\omega$ and  $h_{\o}$
are $G$-invariant since $\o_{KE}$ and $\phi$ are, and it follows that the
family $\phi_t$ is $G$-invariant.  Since the Ricci flow preserves
the $G$-invariance, $u_t$ and $\psi_t$ are $G$-invariant.  The proof of Step 1 then goes through
as before.  Step 2 follows as above once we observe that $\phi_t - \phi_1$ is $G$-invariant.  Q.E.D.

\bigskip
\noindent
{\it Acknowledgements:} We would like to thank the referee for pointing out some inaccuracies
in the original version of this manuscript.

\end{document}